\documentclass[a4paper,12pt]{article}
\usepackage[all]{xy}
\usepackage[T1]{fontenc}
\usepackage[dvips]{graphicx}
\usepackage{amsmath,amssymb,amsfonts,amsthm,latexsym,mathrsfs,textcomp,verbatim}
\xyoption{web}

\title{Site characterizations for\\ geometric invariants of toposes}
\author{{Olivia Caramello} \vspace{3 mm}\\ {\small DPMMS, University of Cambridge,}\\{\small Wilberforce Road, Cambridge CB3 0WB, UK}\\{\small O.Caramello@dpmms.cam.ac.uk}}
\date{December 12, 2011}
\begin{document}
\bgroup           % fake a titlepage 
\let\footnoterule\relax  % no rule above thanks footnotes 
\maketitle
\begin{abstract}
We discuss the problem of characterizing the property of a Grothen-\\dieck topos to satisfy a given `geometric' invariant as a property of its sites of definition, and indicate a set of general techniques for establishing such criteria. We then apply our methodologies to specific invariants, notably including the property of a Grothendieck topos to be localic (resp. atomic, locally connected, equivalent to a presheaf topos), obtaining explicit site characterizations for them.
\end{abstract} 
\egroup 
\vspace{5 mm}

%MACROS-----------------------------------------------------------------------------------------------------------------------

%	European dates ``19 April 1990'' not ``April 19, 1990''
\def\Monthnameof#1{\ifcase#1\or
   January\or February\or March\or April\or May\or June\or
   July\or August\or September\or October\or November\or December\fi}
\def\today{\number\day~\Monthnameof\month~\number\year}

%===========================================================================
%	END OF PROOF BOX
%
%
%  The complexity of the macro necessary to get a little box on the
%  right-hand-side at the end of a proof is amazing.  It really does
%  have to be this long!  Otherwise you're liable to get it at the
%  beginning of the next line, or even on the next page.
%
\def\pushright#1{{%        set up
   \parfillskip=0pt            % so \par doesnt push \square to left
   \widowpenalty=10000         % so we dont break the page before \square
   \displaywidowpenalty=10000  % ditto
   \finalhyphendemerits=0      % TeXbook exercise 14.32
  %
  %                 horizontal
   \leavevmode                 % \nobreak means lines not pages
   \unskip                     % remove previous space or glue
   \nobreak                    % don't break lines
   \hfil                       % ragged right if we spill over
   \penalty50                  % discouragement to do so
   \hskip.2em                  % ensure some space
   \null                       % anchor following \hfill
   \hfill                      % push \square to right
   {#1}                        % the end-of-proof mark (or whatever)
  %
  %                   vertical
   \par}}                      % build paragraph

% prefer proofs with statements, also space after
\def\qed{\pushright{$\square$}\penalty-700 \smallskip}

\newtheorem{theorem}{Theorem}[section]

\newtheorem{proposition}[theorem]{Proposition}

\newtheorem{scholium}[theorem]{Scholium}

\newtheorem{lemma}[theorem]{Lemma}

\newtheorem{corollary}[theorem]{Corollary}

\newtheorem{conjecture}[theorem]{Conjecture}

\newenvironment{proofs}%
 {\begin{trivlist}\item[]{\bf Proof }}%
 {\qed\end{trivlist}}

  \newtheorem{rmk}[theorem]{Remark}
\newenvironment{remark}{\begin{rmk}\em}{\end{rmk}}

  \newtheorem{rmks}[theorem]{Remarks}
\newenvironment{remarks}{\begin{rmks}\em}{\end{rmks}}

  \newtheorem{defn}[theorem]{Definition}
\newenvironment{definition}{\begin{defn}\em}{\end{defn}}

  \newtheorem{eg}[theorem]{Example}
\newenvironment{example}{\begin{eg}\em}{\end{eg}}

  \newtheorem{egs}[theorem]{Examples}
\newenvironment{examples}{\begin{egs}\em}{\end{egs}}

%%%%%%%%%%%%%%%%%%%%%%%%%%%%%%%%%%%%%%%%%%%%%%%%%%%%%%%%%%%%%%%%%%%%%%

%  change some single-character symbols to be more appropriate for logic
\mathcode`\<="4268  % < = \langle
\mathcode`\>="5269  % > = \rangle
\mathcode`\.="313A  % make . a binary  relation
\mathchardef\semicolon="603B % the original
\mathchardef\gt="313E
\mathchardef\lt="313C

\newcommand{\biimp}% bi-implication
 {\!\Leftrightarrow\!}

\newcommand{\bjg}% bi-judgement
 {\mathrel{{\dashv}\,{\vdash}}}

\newcommand{\cod}% codomain
 {{\rm cod}}

\newcommand{\dom}% domain
 {{\rm dom}}

\newcommand{\epi}% epimorphism
 {\twoheadrightarrow}

\newcommand{\Ind}[1]% ind-completion of #1
 {{\rm Ind}\hy #1}

\newcommand{\mono}% monomorphism 
 {\rightarrowtail}

\newcommand{\nml}% normal subgroup
 {\triangleleft}

\newcommand{\ob}% class of objects
 {{\rm ob}}

\newcommand{\op}% opposite category
 {^{\rm op}}

\newcommand{\pepi}% partial epimorphism
 {\rightharpoondown\kern-0.9em\rightharpoondown}

\newcommand{\pmap}% partial map arrow
 {\rightharpoondown}

\newcommand{\Set}% category of sets
 {{\bf Set }}

\newcommand{\Sh}% category of sheaves
 {{\bf Sh}}

\newcommand{\sh}% category of sheaves
 {{\bf sh}}

\newcommand{\Sub}% subobject lattice
 {{\rm Sub}}

\section{Introduction}

In \cite{OC10} we advocate that Grothendieck toposes can effectively act as unifying spaces in Mathematics serving as `bridges' for transferring information between distinct mathematical theories. The transfer of information between theories classified by the same topos, represented by different sites of definition $({\cal C}, J)$ and $({\cal C}', J')$ of their classifying topos, takes place by expressing topos-theoretic invariant properties (resp. constructions) on the topos in terms of properties (resp. constructions) of its two different sites of definition. For each invariant, we thus have a `bridge' 
\vspace{0.5cm}

\[  
\xymatrix {
 & & \Sh({\cal C}, J) \simeq \Sh({\cal C}', J') \ar@/^12pt/@{-->}[drr] & & \\
({\cal C}, J)  \ar@/^12pt/@{-->}[urr] & & & & ({\cal C}', J')}
\]
whose central part is the equivalence of toposes $\Sh({\cal C}, J) \simeq \Sh({\cal C}', J')$ and whose `legs', represented by the dashed arrows, are given by the site characterizations corresponding to the given invariant.

In light of this view, it becomes essential to be able to establish, for relevant classes of topos-theoretic invariant properties $I$, criteria of the kind `a topos $\Sh({\cal C}, J)$ satisfies the property $I$ if and only if the site $({\cal C}, J)$ satisfies a property $P_{({\cal C}, J)}$ (explicitly written in the language of the site $({\cal C}, J)$)', holding for any site $({\cal C}, J)$ or for appropriate classes of sites. Indeed, any such criterion gives us the possibility, in presence of any Morita-equivalence, to operate an automatic transfer of information between the two theories, leading to concrete mathematical results of various nature. Particular cases of these general results have already been applied by the author in several different contexts (cf. for example \cite{OC2} and \cite{OC8}), and in fact the primary aim of this paper is to make a systematic investigation of the problem of obtaining site characterizations which can be conveniently applied in connection to our general philosophy `toposes as bridges'.     

Throughout the past years, site characterizations have been established for several important geometric invariants of toposes, including the ones considered in the present paper (cf. in particular \cite{atomic} and \cite{molecular}, and \cite{El} as a general reference); however, all of these characterizations are of form `A Grothendieck topos satisfies an invariant $I$ if and only if \emph{there exists} a site of definition $({\cal C}, J)$ of it satisfying a certain property $P_{({\cal C}, J)}$'; as such, they cannot be directly applied in connection to the philosophy `toposes as bridges', since the criteria that they give rise to only allow one to enter a given bridge (i.e., to pass from the property $P_{({\cal C}, J)}$ of the site $({\cal C}, J)$ to the invariant $I$) and not to exit from it. In fact, not even the proofs of these results provide information which one can exploit to obtain site characterizations going in the other direction. Therefore, to achieve our goal, the problem needs to be completely reconsidered and approached from a different angle; we do so in this paper, by adopting the point of view of separating sets of toposes. In fact, it turns out that most of the geometric invariants of toposes considered in the literature, notably including the property of a topos to be localic (resp. atomic, locally connected, equivalent to a presheaf topos, coherent), can be expressed in terms of the existence of a separating set of objects for the topos satisfying some invariant property. In this paper, we show that expressing topos-theoretic invariants in terms of the existence of a separating set of objects of the topos satisfying some property paves the way for natural `unravelings' of such invariants as properties of the sites of definition of the topos, and hence for criteria of the desired form.

The paper is organized as follows. In section \ref{geometry}, we make some general remarks about the problem of obtaining site characterizations for geometric invariants of toposes, while in the following sections we apply these considerations to the specific invariants mentioned above, obtaining natural site characterizations for them of the desired kind. Besides their technical interest, these results are meant to provide the reader with a general idea of how the technique `toposes as bridges' introduced in \cite{OC10} actually works in a variety of different cases.    

This work should be considered as a companion to \cite{OC8}, where several syntactic characterizations of geometric invariants on toposes in terms of the theories classified by them were obtained.       
    
\subsection{Terminology and notation}

Our terminology and notation is borrowed from \cite{El}, if not otherwise indicated. 

Moreover, we will employ the following conventions. 

Given a Grothendieck site $({\cal C}, J)$, we denote by $a_{J}:[{\cal C}^{\textrm{op}}, \Set]\to \Sh({\cal C}, J)$ the associated sheaf functor, and by $\eta$ the unit of the adjunction between $a_{J}$ and the canonical inclusion $\Sh({\cal C}, J) \hookrightarrow [{\cal C}^{\textrm{op}}, \Set]$. We denote by $l:{\cal C}\to \Sh({\cal C}, J)$ the composite of the Yoneda embedding $Y:{\cal C}\to [{\cal C}^{\textrm{op}}, \Set]$ with the associated sheaf functor $a_{J}:[{\cal C}^{\textrm{op}}, \Set]\to \Sh({\cal C}, J)$.

To mean that $c$ is an object of a category $\cal C$, we simply write $c\in {\cal C}$.

All the toposes considered in this paper will be Grothendieck toposes, if not otherwise stated.

\section{Geometric invariants of toposes}\label{geometry}

Several topologically-inspired invariants of Grothendieck toposes have been considered in the literature. In fact, as emphasized by Grothendieck himself, a topos can be conveniently considered as a generalized space apt to be studied by adopting a topological intuition. Indeed, a topos $\Sh({\cal C}, J)$ can be seen as a sort of completion of the site $({\cal C}, J)$, on which one can define invariants which correspond to (in the sense of being logically equivalent, or implied by) natural `geometric' properties of sites, thus representing analogues, in the topos-theoretic setting, of classical properties of topological spaces. 

The natural topos-theoretic analogue of the notion of basis of a topological space is the notion of separating set of objects of a topos. Indeed, a basis of a topological space $X$ can be considered as a full subcategory $\cal B$ of the the poset category ${\cal O}(X)$ of open sets of $X$ which is $J_{X}$-dense, where $J_{X}$ is the canonical topology on ${\cal O}(X)$; similarly, a separating set of a topos $\cal E$ can be regarded as a full subcategory of $\cal E$ which is $J_{\cal E}$-dense, where $J_{\cal E}$ is the canonical topology on the topos $\cal E$. Note that if $\cal B$ is a basis of a topological space $X$ then
\[
\Sh(X)\simeq \Sh({\cal B}, J_{X}|_{\cal B}),
\]
by Grothendieck's Comparison Lemma; similarly, if $\cal C$ is a separating set of a topos $\cal E$ then 
\[
{\cal E}\simeq \Sh({\cal C}, J_{\cal E}|_{\cal C}).
\]       
The notion of separating set is intimately related to that of site; in fact, the separating sets of a topos $\cal E$ are (up to the obvious notion of isomorphism) precisely the sets of the form $L_{({\cal C}, J)}:=\{l(c) \textrm{ | } c\in {\cal C}\}$, where $({\cal C}, J)$ is a site of definition of $\cal E$ and $l:{\cal C} \to \Sh({\cal  C}, J)$ is the functor given by the composite of the Yoneda embedding with the associated sheaf functor. 

We can thus naturally expect the properties of (sober) topological spaces which can be expressed in terms of the existence of a basis for the space satisfying a certain property $P$ to be naturally generalizable to the topos-theoretic setting, by replacing bases with separating sets and the property $P$ with an appropriate topos-theoretic analogue. For example, the property of a topos to be atomic (resp. locally connected, coherent) represents a natural topos-theoretic analogue of the property of a space to be discrete (resp. locally connected, coherent). 

Of course, not all the natural properties of topological spaces can be expressed in terms of the existence of a basis satisfying a certain condition; many can be expressed as frame-theoretic properties $P$ of the top element of their frame of open sets. In these cases, natural topos-theoretic analogues of them can be obtained by replacing the top element of the frame of open sets with the terminal object of the topos and the property $P$ with an appropriate topos-theoretic analogue. For example, the property of a topos to be two-valued (resp. compact) represents a natural topos-theoretic analogue of the property of a topological space to be trivial (resp. compact).   

In this paper we shall only be concerned with invariants of the first kind, that is of the form `to have a separating set of objects satisfying some property $P$', but the techniques that we shall elaborate will be also adaptable to invariants of the second kind. 

The problem of finding effective site characterizations for invariants of one kind or another of course admits a satisfactory solution or not depending on the specific invariant under consideration; nonetheless, one can identify some properties which are responsible for such invariants to admit explicit site characterizations. 

A general remark which, as we will see, turns out to be extremely useful in practice is the following: if the property $P$ descends along epimorphisms (that is, for any epimorphism $f:A\to B$ in the topos, if $A$ satisfies $P$ then $B$ satisfies $P$) then one can try to obtain an explicit site characterization of the invariant `to have a separating set of objects satisfying the property $P$', as follows. A topos $\Sh({\cal C}, J)$ has a separating set of objects satisfying $P$ if and only if every object of $\Sh({\cal C}, J)$ of the form $l(c)$ is covered by an epimorphic family of arrows whose domains satisfy $P$. Clearly, if the property $P$ descends along epimorphisms then we can suppose, without loss of generality, this family of arrows to consist entirely of monomorphisms, which can be supposed to be, up to isomorphism, of the form $a_{J}(S\mono {\cal C}(c,-))$ for some ($J$-closed) sieve $S$ on $c$. Therefore, provided that the property $P$ is sufficiently well-behaved to the extent of admitting an `unraveling' of the condition of an object of the kind $a_{J}(S)$ to satisfy $P$ as an explicit condition on $S$ written in the language of the site $({\cal C}, J)$, we have an explicit site characterization for our invariant of the required form. Examples of application of this method are given in sections \ref{localic}, \ref{atomic}, \ref{locallyconnected} and \ref{other}. 
    
Of course, if the property $P$ does not descend along epimorphisms, it still make sense to look for explicit site characterizations for the given invariant, holding for large classes of sites if not for all the sites; an example is given by the property of a topos to be equivalent to a presheaf topos, which is expressible as the requirement `to have a separating set of irreducible objects', which we shall discuss in section \ref{presheaf}. Anyway, one should be aware that the process of obtaining characterizations for these invariants is much less canonical than the one described above for the invariants corresponding to properties $P$ which descend along epimorphisms.

It must also be said that many geometrically-motivated invariant properties of toposes can be naturally expressed in terms of the existence of a separating set for the topos satisfying a property which cannot be expressed as the requirement that all the objects in the separating set satisfy a certain condition. For example, the property of a topos to be coherent can be expressed as the existence of a separating set of compact objects which is closed under finite limits in the topos (cf. \cite{OC8}). Invariants of this kind are in general more hardly tractable, from the point of view of site characterizations, than those discussed above; nonetheless, as witnessed by the past literature on the subject (cf. for example \cite{El}), partial or even complete characterizations for them (holding for large classes of sites) can be achieved, by exploiting the particular `combinatorics' of the sites of definition of the topos in relation to the given invariant.

\section{Localic toposes}\label{localic}

In this section we shall address the problem of finding bijective site characterizations for the property of a topos to be localic.

Recall that a Grothendieck topos $\cal E$ is said to be \emph{localic} if it has a separating set of subterminal objects. We seek criteria for a topos $\Sh({\cal C}, J)$ of sheaves on a site $({\cal C}, J)$ to be localic. 

We start by observing that the objects of the form $a_{J}({\cal C}(-,c))$ are a separating set for $\Sh({\cal C}, J)$; so this topos is localic if and only if for every $c\in {\cal C}$ the family of arrows from subterminal objects to $a_{J}({\cal C}(-,c))$ is jointly epimorphic. 
Now, the subterminal objects of the topos $\Sh({\cal C}, J)$ can be identified with the $J$-ideals on $\cal C$. An arrow $I \to a_{J}({\cal C}(-,c))$ in $\Sh({\cal C}, J)$, where $c$ is an object of $\cal C$ and $I$ is a $J$-ideal on $\cal C$, is the image under $a_{J}$ of an arrow $I\to {\cal C}(-,c)$ in $[{\cal C}^{\textrm{op}}, \Set]$. Such an arrow can be described concretely as a function which assigns to every element $d\in I$ an arrow $\alpha(d):d\to c$ in such a way that for any arrow $g:d\to d'$ in $\cal C$ between elements of $I$, $\alpha(d')\circ g=\alpha(d)$; we shall refer to such an arrow as a $I$-matching family on $c$. We also observe that the cover-mono factorization of an arrow $\alpha:I\to {\cal C}(-,c)$ in the topos $[{\cal C}^{\textrm{op}}, \Set]$ yields a subobject of ${\cal C}(-,c)$, in other words a sieve $S_{\alpha}$ on $c$, defined by the formula $S_{I}=\{\alpha(d) \textrm{ | } d\in I\}$ (note that this is a sieve since $\alpha$ is a matching family). Notice that, since $I$ is a subterminal object in $\Sh({\cal C}, J)$, every arrow from $I$ to $a_{J}({\cal C}(-,c))$ is monic, and hence for any $I$-matching family $\alpha$ on $c$, $a_{J}(\alpha)$ is isomorphic to $a_{J}(S_{I}\mono {\cal C}(-,c))$. 

Therefore, we obtain the following criterion for $\Sh({\cal C}, J)$ to be localic.

\begin{theorem}\label{critlocalic}
Let $({\cal C}, J)$ be a site. Then, with the above notation, the topos $\Sh({\cal C}, J)$ is localic if and only if for every $c\in {\cal C}$ there exists a family ${\cal F}_{c}$ of $J$-ideals on ${\cal C}$ and for every ideal $I\in {\cal F}_{c}$ a $I$-matching family $\alpha_{I}$ on $c$ such that the sieve $\{\alpha_{I}(d):d\to c \textrm{ for some } d\in I \textrm{ and } I \in {\cal F}_{c} \}$ is $J$-covering.
\end{theorem}\qed    

It is interesting to apply the theorem to presheaf toposes and to toposes of sheaves on a geometric site.

\begin{corollary}
Let $\cal C$ be a small category. Then the topos $[{\cal C}^{\textrm{op}}, \Set]$ is localic if and only if $\cal C$ is a preorder.
\end{corollary}

\begin{proofs}
We can apply Theorem \ref{critlocalic} by regarding $[{\cal C}^{\textrm{op}}, \Set]$ as the topos $\Sh({\cal C}, J)$ where $J$ is the trivial topology on $\cal C$. Notice that the $J$-ideals on $\cal C$ are in this case simply the ideals on $\cal C$, that is the sets $I$ of objects of $\cal C$ such that for any arrow $f:a\to b$ in $\cal C$, if $b\in I$ then $a\in I$. The condition of the criterion says that for every $c\in {\cal C}$ there is an ideal $I$ such that $c\in I$ and a $I$-matching family $\alpha_{I}$ on $c$ such that $\alpha_{I}(c)=1_{c}$. Notice that the latter condition implies that for any $d\in I$ and any arrow $g:d\to c$ in $\cal C$, $\alpha_{I}(d)=\alpha_{I}(c)\circ g=1_{c}\circ g=g$; and from this calculation it is clear that a $I$-matching family on $c$ exists if and only if for every element $d\in I$ there is exactly one arrow $d\in c$ in $\cal C$. Our thesis thus immediately follows. 
\end{proofs}

\begin{corollary}
Let $({\cal C}, J_{\cal C})$ be a geometric site. Then the topos $\Sh({\cal C}, J_{\cal C})$ is localic if and only if for every object $c\in {\cal C}$ there exists a covering family of arrows $\{f_{i}:dom(f_{i})\to c \textrm{ | } i\in I \}$ such that for every $i\in I$ and every object $d\in {\cal C}$ and every arrows $g,h:d\to dom(f_{i})$, $f_{i}\circ g=f_{i}\circ h$. 
\end{corollary}

\begin{proofs}
It suffices to notice, using Theorem \ref{critlocalic}, that, by definition of geometric topology, the $J$-ideals are precisely those generated by a single object, and hence that a $I$-matching family on an object $c$ for such an ideal $I$ generated by an object $a$ of $\cal C$ corresponds to an arrow $f:a\to c$ with the property that for any object $b$ and any arrows $g,h:b\to a$ in $\cal C$, we have $f\circ g=f\circ h$.
\end{proofs}

\section{Atomic toposes}\label{atomic}

In this section we shall investigate atomic toposes from the point of view of their site characterizations.

Given a topos $\cal E$, we recall that an object $a$ of $\cal E$ is said to be an \emph{atom} of $\cal E$ if the only subobjects of $a$ in $\cal E$ are the identity subobject and the zero one, and they are distinct from each other.

The following lemma, expressing the fact that atoms descend along epimorphisms, will be useful to us.

\begin{lemma}\label{epiatom}
Let $f:a\to b$ be an epimorphism in a topos $\cal E$. If $a$ is an atom then $b$ is an atom.
\end{lemma}

\begin{proofs}
Let $m:b'\to b$ be a monomorphism in $\cal E$. Consider the pullback in $\cal E$ of $m$ along the arrow $f$.
\[  
\xymatrix {
f^{\ast}(b') \ar[d]_{f^{\ast}(m)} \ar[r]^{g} & b'  \ar[d]^{m} \\
a \ar[r]_{f} & b } 
\] 
The arrow $f^{\ast}(m):f^{\ast}(b')\to a$ is a monomorphism (being the pullback of a monomorphism), and hence either $f^{\ast}(m)=1_{a}$ or $f^{\ast}(m)$ is equal to the zero subobject of $a$. On the other hand, the arrow $g$ is an epimorphism, it being the pullback of an epimorphism. So, if $f^{\ast}(m)=1_{a}$ then by the uniqueness (up to isomorphism) of the epi-mono factorizations of arrows in a topos it follows that $m$ is an isomorphism, while if $f^{\ast}(b')\cong 0$ then we have an epimorphism $g:0\to b'$ and hence $b'\cong 0$.  
\end{proofs}

We shall also need the following proposition. 

\begin{proposition}\label{propatom}
Let $({\cal C}, J)$ be a site and $S$ be a $J$-closed sieve on an object $c$ of $\cal C$. Then $a_{J}(S)$ is an atom of $\Sh({\cal C}, J)$ if and only if $\emptyset \notin J(dom(f))$ for some $f\in S$, and for every subsieve $S'\subseteq S$, either for every $f\in S$ $f^{\ast}(S')\in J(dom(f))$ or for every $g\in S'$ $\emptyset \in J(dom(g))$.

In particular, $l(c)$ is an atom if and only if for every sieve $S$ on $c$, either for every arrow $f\in S$, $f^{\ast}(S)\in J(dom(f))$ or for every arrow $g\in S$, $\emptyset \in J(dom(g))$. 
\end{proposition}

\begin{proofs}
We start observing that for any local operator $j$ on a topos $\cal E$, with associated sheaf functor $a_{j}:{\cal E}\to \sh_{j}({\cal E})$, an object $a$ of $\cal E$ satisfies $a_{j}(a)\cong 0_{\sh_{j}({\cal E})}$ if and only if, denoted by $a\epi a' \mono 1$ the epi-mono factorization of the unique arrow $a\to 1$ in $\cal E$, $a_{j}(a')\cong 0_{\sh_{j}({\cal E})}$; this immediately follows from the fact that for any epimorphism in a topos its domain is isomorphic to zero if and only if its codomain is isomorphic to zero, combined with the observation that the associated sheaf functor preserves epimorphisms. Note in passing that this remark can be used to obtain explicit characterizations of the presheaves which are sent to the zero sheaf by a given associated sheaf functor.

Given a site $({\cal C}, J)$, we would like to understand when a certain sieve $S$ on an object $c\in {\cal C}$ has the property that $a_{J}(S)\cong 0$. We consider the epi-mono factorization $S\epi I_{S}\mono 1$ of the unique arrow $S\mono {\cal C}(-,c)\to 1$ in $[{\cal C}^{\textrm{op}}, \Set]$. Since $1$ is a $J$-sheaf the associated sheaf functor applied to a subterminal object $U$ coincides with its $J$-closure $c_{J}(U)$. Therefore $a_{J}(S)\cong c_{J}(I_{S})$. Now, $I_{S}$ is clearly given by the formula $I_{S}:\{dom(f) \textrm{ | } f\in S \}$, from which it follows, recalling that the zero subterminal in $\Sh({\cal C}, J)$ corresponds to the ideal $\{c\in {\cal C} \textrm{ | } \emptyset \in J(c)\}$, that $a_{J}(S)\cong 0$ if and only if for every $f\in S$, $\emptyset\in J(dom(f))$. 

Now we want to investigate under what conditions $a_{J}(S)$ is an atom of $\Sh({\cal C}, J)$. We have already characterized the conditions that make it non-zero. We observe that every subobject in $\Sh({\cal C}, J)$ of $l(c)$ is of the form $a_{J}(i)$ for some inclusion of subsieves $i:S'\subseteq S$, and we have that $a_{J}(i)$ is an isomorphism if and only if $i$ is $c_{J}$-dense (as a subobject in $[{\cal C}^{\textrm{op}},\Set]$), equivalently if $c_{J}(S)'\cong c_{J}(S)$ (or, alternatively, $S\subseteq c_{J}(S')$). The thesis thus follows immediately from the explicit description of the closure operator $c_{J}$.  
\end{proofs}

\begin{remark}\label{rematom}
Notice that if $a_{J}(S)$ is an atom of $\Sh({\cal C}, J)$ and $f\in S$ is such that $\emptyset \notin J(dom(f))$ then $a_{J}((f))\ncong 0$ and hence $a_{J}((f))\cong a_{J}(S)$. So, from the proof of the proposition we see that, to check that $a_{J}(S)$ is an atom it is equivalent to verify that for any arrow $g$ which factors through $f$, either the sieve generated by it is sent by $a_{J}$ to zero (equivalently, $\emptyset \in J(dom(g))$) or $f^{\ast}((g))\in J(dom(f))$. 
\end{remark}

Now we would like to understand in concrete terms what it means for a topos $\Sh({\cal C}, J)$ to have a separating set of atoms. First, we observe that this condition is equivalent to saying that for every $c\in {\cal C}$, $l(c)$ can be covered by an epimorphic family of arrows whose domains are atoms of $\Sh({\cal C}, J)$. By Lemma \ref{epiatom}, we can suppose without loss of generality (by possibly replacing any arrow in the given epimorphic family by its image), that all the arrows in the family are monic. Therefore, they are all of the form $a_{J}(S\mono {\cal C}(-, c))$ for some sieve $S$ on $c$. By Remark \ref{rematom}, we can suppose that $S$ to be the $J$-closure of a sieve generated by a single arrow. Therefore, using Proposition \ref{propatom}, we obtain the following characterization theorem.

\begin{theorem}
Let $({\cal C}, J)$ be a site. Then the topos $\Sh({\cal C}, J)$ is atomic if and only if for every $c\in {\cal C}$ there exists a $J$-covering sieve on $c$ generated by arrows $f$ with the property that $\emptyset \notin J(dom(f))$ and for every arrow $g$ which factors through $f$, either $\emptyset \in J(dom(g))$ or $f^{\ast}((g))\in J(dom(f))$. 
\end{theorem}\qed

\section{Locally connected toposes}\label{locallyconnected}

In this section we address the problem of establishing site characterizations for locally connected toposes.

Recall that a Grothendieck topos is locally connected if the inverse image functor of the unique geometric morphism from the topos to $\Set$ has a left adjoint. By the results in \cite{OC8}, locally connected Grothendieck toposes can be equivalently characterized as the Grothendieck toposes which have a separating set of indecomposable objects. Recall that a object of a Grothendieck topos is said to be indecomposable if it does not admit any non-trivial (set-indexed) coproduct decompositions.
 
Let us start with a lemma, which expresses the fact that indecomposable objects `descend' along epimorphisms.

\begin{lemma}\label{indepic}
Let $f:a\to b$ be an epimorphism in a Grothendieck topos $\cal E$. If $a$ is indecomposable then $b$ is indecomposable.
\end{lemma}

\begin{proofs}
This follows immediately from the fact that if $f:a\to b$ is an epimorphism then the pullback functor ${\cal E}\slash b \to {\cal E}\slash a$ is logical and conservative (cf. \cite{MM}), in light of the fact that coproducts in a topos can be characterized as epimorphic families of pairwise disjoint subobjects.
\end{proofs} 

Let us now turn to the problem of characterizing the indecomposable objects in a general topos $\Sh({\cal C}, J)$.

\begin{proposition}
Let $({\cal C}, J)$ be a site. Then an object of the form $l(S)$, where $S$ is a $J$-closed sieve on an object $c$ of $\cal C$, is indecomposable if and only if the sieve $S$ satisfies the property that for any family $\{S_{i} \textrm{ | } i\in I\}$ of subsieves $S_{i}\subseteq S$ such that for any distinct $i, i'$ and any $f\in S_{i}\cap S_{i'}$, $\emptyset\in J(dom(f))$, if the union $S'$ of the $S_{i}$ is $c_{J}$-dense in $S$ (i.e., for any arrow $f\in S$, $f^{\ast}(\mathbin{\mathop{\textrm{\huge $\cup$}}\limits_{i\in I}}S_{i})\in J(dom(f))$) then some $S_{i}$ is $c_{J}$-dense in $S$ (i.e., for any arrow $f$ in $S$, $f^{\ast}(S_{i})\in J(dom(f))$). 
\end{proposition}
 
\begin{proofs}
It suffices to recall that coproducts in a topos can be characterized as epimorphic families of pairwise disjoint subobjects, and observe that, up to isomorphism, any subobject of $l(S)$ in $\Sh({\cal C}, J)$ is, up to isomorphism, of the form $a_{J}(i):a_{J}(T)\to a_{J}(S)$ where $i$ is the canonical inclusion of a subsieve $T$ of $S$ into $S$.
\end{proofs} 

We shall call a sieve $S$ satisfying the property in the statement of the proposition a $J$-indecomposable sieve.

Using this Proposition, we can easily get a site characterization for a topos $\Sh({\cal C}, J)$ to be locally connected.

\begin{theorem}
Let $({\cal C}, J)$ be a site. Then the topos $\Sh({\cal C}, J)$ is locally connected if and only if for every $c\in {\cal C}$ there exists a family $\{S_{i} \textrm{ | } i\in I\}$ of $J$-closed $J$-indecomposable sieves on $c$ such that the union $\mathbin{\mathop{\textrm{\huge $\cup$}}\limits_{i\in I}}S_{i}$ is $J$-covering. 
\end{theorem}

\begin{proofs}
The topos $\Sh({\cal C}, J)$ is locally connected if and only if it has a separating set of indecomposable objects, equivalently for every $c\in {\cal C}$, the family of arrows from indecomposable objects to $l(c)$ is epimorphic. By Lemma \ref{indepic}, we can suppose without loss of generality that the arrows belonging to this family are monic, and hence are, up to isomorphism, of the form $a_{J}(i):a_{J}(S)\mono l(c)\cong a_{J}({\cal C}(-,c))$ where $i$ is the canonical inclusion $S\mono {\cal C}(-, c)$ of a $J$-closed sieve $S$ into ${\cal C}(-, c)$. Clearly a family of arrows of this form is epimorphic in $\Sh({\cal C}, J)$ if and only if the union of the corresponding sieves is $J$-covering on $c$, from which our thesis follows. 
\end{proofs}

\section{Toposes which are equivalent to a\\ presheaf topos}\label{presheaf}

In this section we consider the invariant property of a topos to be equivalent to a presheaf topos in relation to the problem of obtaining site characterizations for it.

Let us start with a technical lemma.

\begin{lemma}\label{lemmairr}
Let $({\cal C}, J)$ be a site and $S$ be a $J$-closed sieve on an object $c$ of $\cal C$. Then for any sieve $R$ on $l(S)$ in $\Sh({\cal C}, J)$, $R$ is covering if and only if $S$ is equal to the $J$-closure of the sieve $\{f\in S \textrm{ | } l^{S}(f)\in R\}$ (where $l^{S}(f)$ denotes the factorization of $l(f)$ through the canonical monomorphism $a_{J}(S)\mono l(c)$). In particular, a sieve $R$ is epimorphic on $l(c)$ if and only if the sieve $\tilde{R}_{S}:=\{f\in S \textrm{ | } l(f)\in R\}$ is $J$-covering. 
\end{lemma}

\begin{proofs}
Given a sieve $R$ on $l(S)$, let us consider its pullback in $[{\cal C}^{\textrm{op}}, \Set]$ along the unit $\eta_{c}$ of the reflection corresponding to the subtopos $\Sh({\cal C}, J)\hookrightarrow [{\cal C}^{\textrm{op}}, \Set]$ at the object ${\cal C}(-, c)$. For any $f\in R$, we can cover $\eta_{c}^{\ast}(l(dom(f)))$ in $[{\cal C}^{\textrm{op}}, \Set]$ with an epimorphic family whose domains are representable functors. By the fullness of the Yoneda Lemma, the arrow obtained by composing any such arrow ${\cal C}(-, d)\to dom(f)$ first with $\eta_{c}^{\ast}(l^{S}(f))$ and then with the canonical monomorphism $S\mono {\cal C}(-, c)$ is of the form ${\cal C}(-, h)$ for some arrow $h:d\to c$ in $\cal C$. Notice that the arrows of the form ${\cal C}(-, h)$ corresponding to a given arrow $f\in R$ are jointly epimorphic on $\eta_{c}^{\ast}(l(dom(f)))$ and hence their images under the associated sheaf functor $a_{J}$ yield a jointly epimorphic family $T_{f}$ on $l(dom(f))$.  

\[  
\xymatrix {
\eta_{c}^{\ast}(l(dom(f))) \ar[rr]^{\eta_{c}^{\ast}(l^{S}(f))} \ar[d] & & S \ar[d] \ar[r] & {\cal C}(-, c) \ar[d]^{\eta_{c}} \\
l(dom(f)) \ar[rr]^{l^{S}(f)} & & l(S) \ar[r] &  l(c) } 
\] 

Now, clearly, $R$ is epimorphic if and only if the sieve $R^{m}$ on $l(S)$ obtained by `multicomposing' $R$ with the sieves $T_{f}$ (for $f\in R$) is epimorphic. Notice that $R^{m}$ is generated by the factorizations through the canonical monomorphism $a_{J}(S)\mono l(c)$ of arrows of the form $l(h)$ where $h$ is an arrow in $\cal C$ with codomain $c$. Define the sieve $A$ on $c$ as 
\[
A:=\{k\in S \textrm{ | } l^{S}(k)\in R\}.
\]
Clearly, the sieve $A_{l}$ of arrows of the form $l^{S}(k)$ for $k\in A$ is contained in $R$ and contains $R^{m}$, from which it follows that it is jointly epimorphic if and only if $R$ (equivalently, $R^{m}$) is. From this our thesis clearly follows, since $A_{l}$ is jointly epimorphic if and only if the image of the canonical monomorphism $A\mono S$ under $a_{J}$ is an epimorphism (equivalently, an isomorphism), that is if and only if $A\mono S$ is $c_{J}$-dense, where $c_{J}$ is the closure operator associated to the Grothendieck topology $J$. 

\end{proofs}

We notice that the particular case of Lemma \ref{lemmairr} when the sieve $S$ is maximal was observed, but not proved, at p. 911 of \cite{El}.

Another useful remark concerning the relationship between sieves and their images under $l$ functors is provided by the following proposition.

\begin{proposition}\label{overline}
Let $({\cal C}, J)$ be a site and $S$ be a sieve on an object $c$ of $\cal C$. Then the sieve 
\[
\overline{S}=\{f:dom(f)\to c \textrm{ | } l(f) \textrm{ factors through } l(g) \textrm{ for some } g\in S \}
\]
is contained in the $J$-closure of $S$. 
\end{proposition}

\begin{proofs}
Suppose that $l(f)$ factors through $l(g)$ for some $g\in S$. The canonical monomorphism $p:f^{\ast}(R)\mono {\cal C}(-, dom(f))$ can be identified with the pullback of the canonical monomorphism $q:S\mono {\cal C}(-, c)$ along ${\cal C}(-, f)$. Hence, if $l(f)$ factors through $l(g)$, $a_{J}(p)$ is isomorphic to the pullback of $a_{J}(q)$ along any factorization of $l(f)$ through $l(g)$; therefore, if $g\in S$ this latter pullback is an isomorphism and hence $p$ is $c_{J}$-dense, in other words $f^{\ast}(R)$ is $J$-covering, i.e. $f$ belongs to the $J$-closure of $S$.   
\end{proofs}

We say that a sieve $S$ on an object $c$ is $l$-closed if $\overline{S}=S$. Note that the $l$-closed sieves are precisely those of the form $\tilde{T}$ for some sieve $T$ on $c$. 

Notice that if $J$ is subcanonical then every sieve is $l$-closed, since the functor $l$ is full and faithful.

We can use the lemma to characterize the irreducible objects of a topos $\Sh({\cal C}, J)$. Recall that an object of a Grothendieck topos is said to be irreducible if every covering sieve on it is maximal. 

Suppose that $P$ is an irreducible object of $\Sh({\cal C}, J)$. We can cover $P$ with a family of arrows whose domains are of the form $l(c)$ for $c\in {\cal C}$; there is thus an arrow of the family, say $e:l(c)\to P$, which is split epic, that is such that there is a monic arrow $m:P\to l(c)$ such that $e\circ m=1_{P}$. Now, $m$ being mono, $P$ is, up to isomorphism, of the form $l(S)$ where $S$ is a sieve on $c$. Now, consider the family of monomorphisms $\{l((f))\mono l(S) \textrm{ | } f\in S\}$. This family covers $l(S)$ and hence, by the irreducibility of $P$, $l(S)\cong l((f))$ for some $f\in S$. Notice that the canonical monomorphism $(f)\mono {\cal C}(-, c)$ is the monic part of the epi-mono factorization of the arrow ${\cal C}(-, f):{\cal C}(-, dom(f))\to {\cal C}(-, c)$ in $[{\cal C}^{\textrm{op}}, \Set]$ and hence $m:l((f))\mono l(c)$ is the monic part of the epi-mono factorization of $l(f):l(dom(f))\to l(c)$ in $\Sh({\cal C}, J)$. Let us denote by $e':l(dom(f))\to l((f))$ the epic part of this factorization (notice that this arrow is in fact split epic by the irreducibility of $l((f))$).

Now, given a sieve $R$ on $l((f))$, $R$ is maximal if and only if $f\in \tilde{R}$. Therefore the condition that for any sieve $R$ on $l((f))$, $R$ is epimorphic if and only if it is maximal can be equivalently expressed as the condition that $\tilde{R}$ is dense in the $c_{J}$-closure of $(f)$ if and only if it is maximal. 

In order to obtain a criterion for irreducibility which does not involve constructions within the topos $\Sh({\cal C}, J)$, we would like to replace the quantification over the sieves $R$ in the topos $\Sh({\cal C}, J)$ with a quantification over sieves in the category $\cal C$. To this end, we investigate whether if $l((f))$ is irreducible then it is the case that for \emph{every} subsieve $S$ of $(f)$ - not just those of the form $\tilde{R}$ for a sieve $R$ on $l((f))$ - $S$ is $c_{J}$-dense in the $c_{J}$-closure of $(f)$ if and only if it is equal to $((f))$, and if not, whether it is possible to characterize intrinsically a class of subsieves which enjoy this property. Notice that $S$ is $c_{J}$-dense in the $c_{J}$-closure of $(f)$ if and only if the sieve $S_{p}$ on $l((f))$ generated by the arrows of the form $l^{(f)}(g)$ for $g\in S$ is covering on $l((f))$, while, by Proposition \ref{overline}, $S_{p}$ is maximal if and if and only if $\overline{S}=(f)$. Therefore the required property is satisfied by all the $l$-closed sieves $S$; that is for any such sieve $S$, $S$ is $c_{J}$-dense in the $c_{J}$-closure of $(f)$ (equivalently, $f^{\ast}(S)\in J(dom(f))$) if and only if it is equal to $(f)$ (equivalently, $f\in S$).

Therefore we can conclude the following result.

\begin{proposition}\label{charirr}
Let $({\cal C}, J)$ be a site and $f$ be an arrow of $\cal C$. Then the object $l((f))$ is irreducible in $\Sh({\cal C}, J)$ if and only if for every $l$-closed sieve $S\subseteq (f)$ on $c$, $f^{\ast}(S)\in J(dom(f))$) if and only if $f\in S$. The requirement of $S$ to be $l$-closed can be omitted if the topology $J$ is subcanonical.      
\end{proposition}\qed 

We now proceed to obtain, by using Proposition \ref{charirr}, an intrinsic site characterization of the toposes which are equivalent to presheaf toposes.

\begin{theorem}\label{irrchar}
Let $({\cal C}, J)$ be a site. Then the topos $\Sh({\cal C}, J)$ has a separating set of irreducible objects (equivalently, is equivalent to a presheaf topos) if and only if for every object $c$ of $\cal C$ there exists a family $\{(k_{i}, w_{i}) \textrm{ | } i\in I\}$ of pairs of composable arrows such that the sieve generated by the family $\{w_{i}\circ k_{i} \textrm{ | } i\in I\}$ is $J$-covering and for every $l$-closed sieve $S\subseteq (k_{i})$ on $cod(k_{i})$, $k_{i}^{\ast}(S)\in J(dom(k_{i}))$ implies $k_{i}\in S$. The requirements of the sieves to be $l$-closed can be omitted if the topology $J$ is subcanonical.       
\end{theorem}

\begin{proofs}
The topos $\Sh({\cal C}, J)$ has a separating set of irreducible objects if and only if every object of the form $l(c)$ for an object $c$ of $\cal C$ is covered by a sieve generated by arrows whose domains are irreducible objects. Given any such arrow $u:P\to l(c)$ in the family (where $P$ is an irreducible object), let us first show that there is a split epic arrow to $P$ of the form $e:l(d)\to P$, with splitting $m:P\mono l(d)$, which, composed with $u$, gives an arrow of the form $l(w)$. Consider the pullback $u':P'\to {\cal C}(-, c)$ of $u$ in $[{\cal C}^{\textrm{op}}, \Set]$ along the unit $\eta_{c}:{\cal C}(-, c)\to l(c)$; $P'$, regarded as an object of $[{\cal C}^{\textrm{op}}, \Set]$, can be covered by an epimorphic family whose domains are representable functors. By the fullness of the Yoneda Lemma, any arrow obtained by composing such an arrow ${\cal C}(-, d)\to P'$ with $u'$ is of the form ${\cal C}(-, w)$ for some arrow $w:d\to c$ in $\cal C$. Notice that the arrows of the form ${\cal C}(-, w)$ are jointly epimorphic on $P'$ and hence their images under the associated sheaf functor $a_{J}$ yield a jointly epimorphic family on $P$, which, by the irreducibility of $P$, must contain a split epic arrow, with monic splitting $m:P\to l(d)$. By the argument preceding Proposition \ref{charirr}, the monomorphism $m:P\to l(d)$ can be identified with the monic part of the epi-mono factorization $l(dom(k))\epi l((k))=P\mono l(d)=l(cod(k))$ of an arrow of the form $l(k)$ where $k$ is an arrow $dom(k)\to d$ in $\cal C$, and $P\cong l((k))$. Let us denote by $z$ the epic part $l(dom(k))\epi l((k))$ of this factorization. Then the arrow $u\circ z$ is equal to $l(w\circ k)$, since $u=u\circ 1_{P}=u\circ e \circ m\circ z=l(w)\circ m \circ z=l(w)\circ l(k)=l(w \circ k)$. Therefore, since the arrow $z$ is epic, the family of the arrows of the form $l(w\circ k)$, where $w$ and $k$ vary as $u$ does in the original epimorphic family, is epimorphic on $l(c)$, equivalently the sieve generated by the family of arrows $\{w\circ k\}$ is $J$-covering on $c$. Thus we can conclude that $\Sh({\cal C}, J)$ has a separating set of irreducible objects if and only if for every object $c$ of $\cal C$ there exists a family $\{(k_{i}, w_{i}) \textrm{ | } i\in I\}$ of pairs of composable arrows such that the sieve generated by the family $\{w_{i}\circ k_{i} \textrm{ | } i\in I\}$ is $J$-covering and $l((k_{i}))$ is irreducible for each $i\in I$. Our thesis now follows by invoking Proposition \ref{charirr}.       
\end{proofs}

\begin{remark}
If $\cal C$ is Cauchy-complete and $J$ is subcanonical then the criterion of Theorem \ref{irrchar} significantly simplifies, as shown in \cite{OC8}, since all the irreducible objects in $\Sh({\cal C}, J)$ are of the form $l(c)$ for some object $c\in {\cal C}$.
\end{remark}

\section{Other invariants}\label{other}

In this section we consider other two invariants on Grothendieck toposes, with the purpose of illustrating additional examples of invariants for which natural site characterizations can be achieved. 

The first invariant that we shall investigate is the property of having a separating set of well-supported objects. Recall that an object $A$ of a topos $\cal E$ is said to be well-supported if the unique arrow $A\to 1$ is an epimorphism. Notice that the property of being well-supported descends along any arrow, that is for any arrow $f:A\to B$ in $\cal E$, if $A$ is well-supported then $B$ is well-supported. 

Let $({\cal C}, J)$ be a site. Clearly, the topos $\Sh({\cal C}, J)$ has a separating set of well-supported objects if and only if every object of the form $l(c)$ (for $c\in {\cal C}$) can be covered by a family of arrows whose domains are well-supported objects. By definition of epimorphic family, if the covering family on $l(c)$ is empty then $l(c)$ is isomorphic to zero, while if the family is non-empty then $l(c)$ is well-supported; from this remark we conclude (classically) that $\Sh({\cal C}, J)$ has a separating set of well-supported objects if and only if for every $c\in {\cal C}$, either $l(c)\cong 0$ or $l(c)$ is well-supported. Now, $l(c)$ is well-supported if and only if every object of $\cal C$ is covered by a $J$-covering sieve generated by arrows whose domains are objects which admit an arrow to $c$, while $l(c)\cong 0$ if and only if $\emptyset\in J(c)$. Therefore, we have the following result.

\begin{theorem}
Let $({\cal C}, J)$ be a site. Then the topos $\Sh({\cal C}, J)$ has a separating set of well-supported objects if and only if for every $c\in {\cal C}$, either $\emptyset \in J(c)$ or for every $d\in {\cal C}$ there exists a sieve $S\in J(d)$ such that for every $f\in S$ there exists an arrow $dom(f)\to c$ in $\cal C$.  
\end{theorem}\qed

We can straightforwardly apply this result to presheaf toposes and to toposes of sheaves on a geometric site.

\begin{corollary}
Let $\cal C$ be a small category. Then the topos $[{\cal C}^{\textrm{op}}, \Set]$ has a separating set of well-supported objects if and only if for any objects $c, d\in {\cal C}$ there exists an arrow $c\to d$ in $\cal C$.
\end{corollary}\qed

\begin{corollary}
Let $({\cal C}, J)$ be a geometric site. Then the topos $\Sh({\cal C}, J)$ has a separating set of well-supported objects if and only if for every $c\in {\cal C}$, either $c\cong 0_{{\cal C}}$ or the unique arrow $c\to 1_{\cal C}$ in $\cal C$ is a cover.
\end{corollary}\qed 

Finally, we consider a fundamental topos-theoretic invariant, namely the property of a topos to be coherent, from the point of view of site characterizations. The property of coherence for a topos is more problematic in relationship to bijective site characterization of the kind we seek than the other invariants considered above in the paper. In \cite{OC8} we observed that a topos is coherent if and only if it has a separating set of compact objects which is closed under finite limits. The problem with this characterization is that we cannot suppose the separating set to be closed under quotients and therefore we cannot apply the usual technique of making the relevant property descend along epimorphisms. On the other hand, the weaker invariant property of having a separating set of compact objects clearly admits a site characterization of the required form, since the property of an object of a topos to be compact descends along epimorphisms. Still, it is possible to achieve bijective site characterizations for the property of a topos to be coherent which hold for large classes of sites (cf. for example \cite{beke} for the case of presheaf toposes). For instance, by using the fact that in a presheaf topos all the representable functors are irreducible objects and that any retract of a coherent object in a coherent topos is coherent, one can immediately deduce that if $[{\cal C}^{\textrm{op}}, \Set]$ is coherent then all the representable functors $y(c)$ are coherent objects; in particular, any finite product $y(c_{1})\times \cdots \times y(c_{n})$ of them is compact and the equalizer of any pair of arrows between such objects is compact (note that both conditions can be expressed as genuine properties of the category $\cal C$; for example, the latter can be expressed by saying that for any arrows $f,g:c\to d$ in $\cal C$, the sieve consisting of all the arrows $h$ with codomain $c$ such that $f\circ h=g\circ h$ is generated by a finite family of arrows). Notice in passing that this kind of characterizations can be profitably applied in presence of any Morita-equivalence involving a presheaf topos according to the philosophy `toposes as bridges' of \cite{OC10}; for example, they allow us to see that the syntactic property of a theory of presheaf type to be coherent has semantic consequences at the level of the `geometry' of its category of finitely presentable models (for instance, the characterization provided by Theorem 2.1 \cite{beke} implies that for any theory of presheaf type $\mathbb T$, if $\mathbb T$ is coherent then its category of finitely presentable models has \textrm{fc} finite limits, in the sense of \cite{beke}). 

This last example is just meant to give the reader an idea of the huge amount of concrete mathematical results in distinct fields that can be obtained, in a `uniform' and essentially automatic way, by using site characterizations for geometric invariants of toposes such as the ones that we have established in the present paper in connection with the philosophy `toposes as bridges'; other notable applications of the same general methodology can be found in \cite{OC2}, \cite{OC8} and \cite{OC10}.

\end{document}